\documentclass{svmult}

\usepackage{makeidx}         
\usepackage{graphicx}        
\usepackage{subfigure}
\usepackage{multicol}        
\usepackage[bottom]{footmisc}

\usepackage{psfrag}

\usepackage{amsmath,amssymb}

\usepackage{color}
\usepackage{fancyhdr}

\DeclareMathOperator{\Id}{Id}

\newcommand{\Keff}{\bold K_{\text{eff}}}
\newcommand{\Keffinv}{\Keff^{-1}}
\newcommand{\uavg}{\langle \bold u \rangle}
\newcommand{\pavg}{\langle p \rangle}

\begin{document}

\title*{Numerical simulation of transport in porous media: some problems from micro to macro scale}
\titlerunning{Transport in porous media}
\author{
Quanji Cai\inst{1},
Sheema Kooshapur\inst{2},
Michael Manhart\inst{2},
Ralf-Peter Mundani\inst{1},
Ernst Rank\inst{1},
Andreas Springer\inst{3},
Boris Vexler\inst{3}
}
\authorrunning{Cai et al.}
%
\institute{
Lehrstuhl Computation in Engineering, Technische Universit\"at M\"unchen\\
         80333 M\"unchen, Germany
\texttt{cai@bv.tum.de;mundani@tum.de;ernst.rank@tum.de}
\and
Fachgebiet Hydromechanik, Technische Universit\"at M\"unchen\\
         80333 M\"unchen, Germany
\texttt{michael.manhart@tum.de;s.kooshapur@bv.tum.de}
\and 
Centre for Mathematical Sciences, M1 Technische Universit\"at M\"unchen\\
         85748 Garching b. M\"unchen, Germany
\texttt{springer@ma.tum.de;vexler@ma.tum.de}
}
%
%

\maketitle

\begin{abstract}
This paper deals with simulation of flow and transport in porous media such as transport of groundwater contaminants. We first discuss how macro scale equations are derived and which terms have to be closed by models. The  transport of tracers  is strongly influenced by pore scale velocity structure and large scale inhomogeneities in the permeability field. The velocity structure on the pore scale is investigated by direct numerical simulations of the 3D velocity field in a random sphere pack. The velocity probability density functions are strongly skewed, including some negative velocities. The large probability for very small velocities might be the reason for non-Fickian dispersion in the initial phase of contaminant transport. We present a method to determine large scale distributions of the permeability field from point-wise velocity measurements. The adjoint-based optimisation algorithm delivers fully satisfying agreement between input and estimated permeability fields. Finally numerical methods for convection dominated tracer transports are investigated from a theoretical point of view. It is shown that high order Finite Element Methods can reduce or even eliminate non-physical oscillations in the solution without introducing additional numerical diffusivity.

\end{abstract}

\pagestyle{empty}
\thispagestyle{fancy}
\lhead{}
\chead{}
\rhead{}
\lfoot{\scriptsize This is a pre-print of an article published in Bader~M., Bungartz~HJ., Weinzierl~T.\ (eds) Advanced Computing. Lecture Notes in Computational Science and Engineering, vol 93, 2013. The final authenticated version is available online at: https://doi.org/10.1007/978-3-642-38762-3\_3}
\cfoot{}
\rfoot{}

\flushbottom

\section{Introduction}
\label{sec:manhartetal_introduction}
For a correct description of reactive flow in porous media, the transport of the reactive species needs to be described correctly. As reaction takes place only in contact zones of the species involved, effective reaction rates are dependent on the microscopic concentration fields which can be strongly heterogeneous. The main problem in predicting concentration fields during tracer transport in a natural porous medium, such as soil, arises from the large range of scales involved. They span from the pore (micro) scale to the field (macro) scale, thus reaching from the range of or smaller than a $\mathrm{\mu m}$ to  the $\mathrm{km}$ range. Thus different techniques are used to simulate tracer transport on different scales.

Transport in porous media is governed by three processes, the advective transport by the macroscopic flow field, the molecular diffusion and the mechanical dispersion due to the randomness of the individual streamlines through the pore space. Modelling dispersion on the macro scale has often been done by assuming an effective diffusivity for the tracer (\cite{bear_72,silva_07}). The resulting advection-diffusion equation can be solved by standard discretisation methods (e.g. FE, FV and FD)  or by stochastic (random walk) methods  (\cite{dentz_04,suciu_13}). Classical (FE, FV and FD) methods lack stability  in advection dominated problems of tracer transport in porous media. Due to sharp gradients and front evolving in the solution, classical non-diffusive tend to produce non-physical oscillations. A way to get rid those oscillations is the introduction of numerical diffusion by upwinding. Another way is to stabilise the FE method by a variational multi-scale formulation \cite{juanes_2005,yang_2009}. 

Modelling the mechanical dispersion by an effective diffusivity needs to regard two aspects, the non-Fickian regime in the initial phase and the dependence of the effective diffusivity on the randomness and structure of the porous matrix (e.g. soil). 
Special methods have been proposed to model non-Fickian dispersion in the initial phase by \cite{dentz_04,hassanizadeh_96,levy_03}. Such methods require knowledge of multi-point/multi-time statistics of the tracer and are therefore difficult to handle. On the other hand, the formulation of effective diffusion coefficients in the Fickian regime also requires knowledge on the randomness of the porous matrix. Preferential paths strongly amplify mechanical dispersion because in relatively slow regions, tracers can stay for a long time. Many studies therefore deal with the description of the permeability fields and their impact on tracer transport (e.g. \cite{dentz_04}). In many cases, the parameters are subject to large uncertainties and can, if at all, only be described stochastically.

Recently, interest has grown in methods relying on velocity probability density functions (PDF). Meyer et al. have proposed a joint velocity-concentration
PDF equation which accounts for advective transport and pore-scale dispersion in porous
media and is solved by a particle method which is able to deal with non-Gaussian
distributions of the velocity field (\cite{meyer_04,meyer_10}); Jenny et al.
introduce a new PDF method for obtaining information about tracer and phase transport by
assuming that the multi-point velocity statistics is known (\cite{jenny_06}).
Nowak et al. show the dependence of hydraulic heads and velocities on the variance of log-conductivity using Monte Carlo simulations. They offer insight into the credibility of first-order second moment
methods for evaluating moments of hydraulic heads. They observe a large deviation of the discharge
components from Gaussian distribution and suggest using more accurate methods such as Monte Carlo
if no assumptions on the shape of distributions are justified (\cite{nowak_08}). Deurer et al. \cite{deurer_04}  measured velocity PDFs in sphere packs by magnetic resonance imaging in various sample volumes to investigate longitudinal and transverse dispersion. They observed a strong dependence of the PDFs from sample volume.

In a research initiative on reactive flows in porous media, three different directions have been followed to improve prediction of concentration fields during the simulation of species transport through a porous medium. Our contributions are in the following fields: (i) proper resolution of the gradients of tracers without numerical diffusion on the macro scale (sec. 5) (ii) description of subfilter fluctuations on micro-scale (sec. 3) description of subfilter fluctuations on macro scale (sec. 4).

{The paper is organised as follows.} In the next section, the  equations describing flow in porous media, both on micro- as on macro-scale are discussed. After that, examples are presented that attack some of the problems in solving these equations by numerical methods. First, pore scale simulations using full solution of the Navier-Stokes equations are presented. Then, a method for parameter identification of an inhomogeneous permeability field is presented. Finally, a high order numerical method for transport on the macro-scale (Darcy-scale) is presented and discussed.


\section{Description of flow in porous media from micro to macro scale} 
In this section some basic quantities on flow in porous media are defined. We start from a definition of the flow quantities on micro- and macro-scale as well as a discussion of the relevant equations of flow and tracer transport. The macro-scale equations are obtained by consequent homogenisation of the micro-scale equations over a representative elementary volume (REV). From this homogenisation, unclosed terms arise that have to be modelled adequately. Some problems of modelling and numerical solution of the respective equations are discussed.

We are considering incompressible  flow of a Newtonian fluid and tracer transport through a porous medium. On the micro-scale, i.e. on volumes as large as the individual pores, the flow is governed by the Navier-Stokes equations, the conservation of mass

\begin{equation}
\label{eq:continuity}
\nabla\cdot\bold u = 0
\end{equation}
and the conservation of momentum
\begin{equation}
\label{eq:navier-stokes}
\rho \partial_t\bold u
+\rho \bold u\cdot\nabla \bold u
=
- \nabla p
+\mu \nabla^2\bold u
\end{equation}  
Here, $\bold u$, $p$, $\rho$ and $\mu$ denote the velocity, pressure, density and dynamic viscosity, respectively. The transport of a tracer in the pore space is described by the convection diffusion equation for the concentration $c$ of the tracer
\begin{equation}
\label{eq:convection-diffusion}
\partial_t c
+ \bold u\cdot\nabla c
=
\Gamma \nabla^2 c\,.
\end{equation}  
Here, $\Gamma$ is the molecular diffusivity.

The formalism of volume averaging   \cite{Whitaker:1986} establishes a rigorous way of deriving macroscopic equations from the microscopic ones. If the total control volume, including fluid and solid phase, is denoted by $V$, then a superficial average of a quantity $\psi$ can be defined the following way
\begin{equation}
\label{eq:superficial-average}
\langle \psi\rangle = \frac{1}{{V}}\int_V\psi dx\,.
\end{equation}
The porosity $\epsilon=V_p/V$ is defined to be the ratio of fluid filled volume (pore space $V_p$) divided by the total volume $V$.
By volume-averaging the momentum equation (\ref{eq:navier-stokes}) the well-known Darcy equation can be obtained
\begin{equation}
\label{eq:darcy}
\langle \bold u\rangle = - {\bold K}\nabla \langle p\rangle\,,
\end{equation} 
in which $\bold K$ denotes the permeability tensor. However, when applying the averaging procedure on a larger scale, the definition of an {\it effective} permeability tensor poses problems as it is not a mere averaging of the permeability tensor at smaller scales. This can be seen by integrating equation (\ref{eq:darcy}) over a larger volume which gives

\begin{equation}
\label{eq:darcy_averaged}
\left\langle\langle \bold u\rangle\right\rangle = 
- \left\langle{\bold K}\nabla \langle p\rangle\right\rangle \neq
- \langle{\bold K}\rangle \nabla\langle p\rangle\,.
\end{equation} 
In measurements, often only large scale permeabilities are accessible, treated as effective permeabilities $K_{\text{eff}}$. If small scale variability of the permeability was accessible, the effective permeability can be obtained by up-scaling methods \cite{Durlofsky_1991}
\begin{equation}
\label{eq:k_eff}
- \left\langle{\bold K}\nabla \langle p\rangle\right\rangle =
- \bold K_{\text{eff}} \nabla\langle p\rangle\,.
\end{equation} 
The dispersion on a macro-scale is dependent on the distribution of the permeabilities on the scale of an REV as this determines whether e.g. preferential flow paths can establish.

When homogenising the convection diffusion equation (\ref{eq:convection-diffusion}), a similar problem arises. Averaging over an REV gives
\begin{equation}
\label{eq:convection-diffusion-averaged}
\partial_t\langle c\rangle
+ \langle\bold u\cdot\nabla c\rangle
=
\Gamma \langle\nabla^2\bold c\rangle\,.
\end{equation}
In here, we have to realise that the second term on the left hand side causes problems, as $\langle\bold u\cdot\nabla c\rangle\neq \langle\bold u\rangle\cdot\nabla \langle c\rangle$. The underlying phenomenon is called dispersion. In most cases, it can be modelled by an additional diffusion using an effective dispersion coefficient \cite{bear_72}
\begin{equation}
\label{eq:convection-diffusion-averaged}
\langle\bold u\cdot\nabla c\rangle
=
\langle\bold u\rangle\cdot\nabla \langle c\rangle
+
\Gamma^{\text{disp}} \nabla^2\langle c\rangle\,.
\end{equation}
An effective dispersion is a good and valuable approach for late phases of tracer transport which are characterised by Gaussian tracer plumes \cite{dentz_04}. Using $\Gamma^{\text{eff}}=\Gamma + \Gamma^{\text{disp}}$, equation (\ref{eq:convection-diffusion-averaged}) is then formulated as 
\begin{equation}
\label{eq:convection-diffusion-averaged-2}
\partial_t\langle c\rangle
+ \langle\bold u\rangle\cdot\nabla \langle c\rangle
=
\Gamma^{\text{eff}} \nabla^2\langle  c\rangle\,.
\end{equation}
These late stages are characterised by Fickian dispersion \cite{dentz_04}. For early phases, strongly non-Gaussian tracer plumes and break-through curves are observed. These stages are characterised by non-Fickian dispersion and need special methods for description.

When flow and transport problems on a macro-scale are addressed the corresponding macroscopic parameters have to be modelled adequately, namely the effective permeability $K_{\text{eff}}$ and the effective dispersion coefficient $\Gamma^{\text{eff}}$. Both can not be directly determined from basic principles. Either empirical correlations, experiments or numerical simulations on the micro-scale have to be used to estimate those macro-scale parameters.

In the following, we present some numerical efforts to improve our understanding of macro-scale parameters and processes. The first one addresses the description of dispersion by knowledge of the micro-scale velocity field, the second one deals with the estimation of the effective permeability distributions by macro-scale measurements and the third effort deals with the solution of the convection diffusion equation in convection dominated transport.


\section{Pore scale simulations of the flow through a random sphere pack.} 
The variability of flow paths and velocities in porous media results in a dispersion of a tracer during its transport through a porous medium. Understanding the variability in the flow field is  the key to understand and model dispersion in a rigorous way. The late phases of dispersion can be modelled by Fickian diffusion with an effective dispersion coefficient, see equation (\ref{eq:convection-diffusion-averaged-2}). Early phases, i.e. non-Fickian transport, need special attention as equation (\ref{eq:convection-diffusion-averaged-2}) can not represent non-Fickian behaviour which is often characterised by strongly skew break-through curves. In the following we present an attempt to understand flow variability in the pore space of a random sphere pack by describing the velocity distribution within the pore space.

We investigate the flow field on the pore scale of regular and random sphere packs
by direct numerical simulation. The full Navier-Stokes equations (\ref{eq:continuity}) and (\ref{eq:navier-stokes}) for an
incompressible, Newtonian fluid are solved by a Finite Volume method on a Cartesian
grid \cite{manhart_02c}. The irregular pore space is represented by an Immersed Boundary Method (IBM) to interpolate the no-slip boundary condition on the spheres to the Cartesian mesh \cite{peller_2010,peller_05}. The spheres are represented by a triangular surface grid of triangle size smaller than the grid spacing of the Cartesian grid. 
The time advancement is done by a low-storage third order Runge-Kutta method \cite{williamson_80}. 
This basic solver is well validated in various flow configurations including laminar and turbulent flows (e.g. \cite{breuer_09,hokpunna_2010,peller_2010}). It has been shown that for viscous flow problems a second order convergence with grid refinement is achieved \cite{peller_2010,peller_05}. The sphere pack is generated by a special algorithm that distributes the spheres randomly in space. To achieve a periodic placement of the spheres, we first arranged the spheres on the faces of the domain. The inner part of the domain is then packed with as many spheres as possible. This method unfortunately results in  a rather more porous area between the faces of the domain and the inner region that has to be taken into account in the post-processing.

We apply periodic boundary conditions in all three space dimensions. The flow is
driven by a constant pressure gradient that is applied as a source term in the momentum equation. The simulation is advanced from rest until convergence has
been reached. As the Reynolds numbers are extremely small, the time to reach
convergence is mainly determined by the diffusion time scale within the pore space.

\paragraph{Grid study}
We checked the accuracy of the method by a convergence study of the flow through a regular sphere pack. In order to obtain the porous geometry we placed 23 spheres in a hexagonal packing arrangement and took out the smallest sized box that would fit into this arrangement and would be periodic in all three directions as our domain. We simulated low Reynolds number flow through this domain which was of size $(L_{x},L_{y},L_{z}) = (4, 2\sqrt{3}, 2\sqrt{√3})\,\mathrm{mm}$. The flow was driven by a pressure gradient of $0. 002\,\mathrm{P_{a}/m}$ in the $x$-direction. The Reynolds number of this setup was in the order of $Re= U_{i}D/{\nu}=1\times10^{-5}$. Here, $D$ is a characteristic length scale such as pore size or sphere diameter, and ${\nu}$ is the kinematic viscosity. $U_{i}$ is the intrinsic velocity, which is defined as the mean pore velocity in the porous domain

 \begin{equation}  
U_{i} = \frac{1}{V_{p}}\int_{V_{p}}^{} u(\textbf{x})\,dv ,
\end{equation}

 with $V_{p}$ being the volume of the pore space. The intrinsic velocity is related to the superficial or Darcy velocity $\langle u\rangle$ by $U_{i}=\langle u\rangle /{\epsilon}$, where $\epsilon$ is the porosity.

\begin{figure}
 \begin{center}
   \subfigure[]{\label{fig:pore-velocity-a}\includegraphics[width=.475\textwidth]{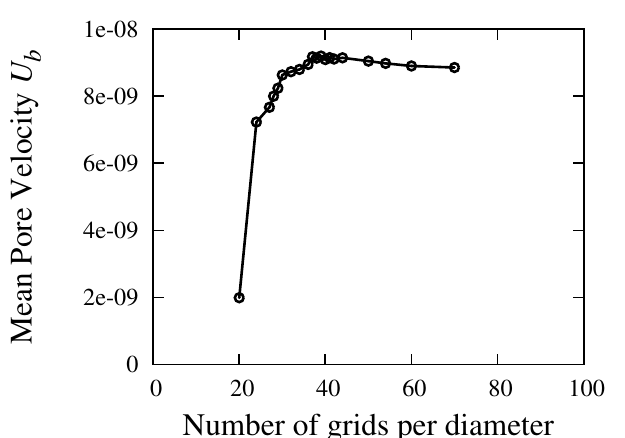}}
   \subfigure[]{\label{fig:pore-velocity-b}\includegraphics[width=.475\textwidth]{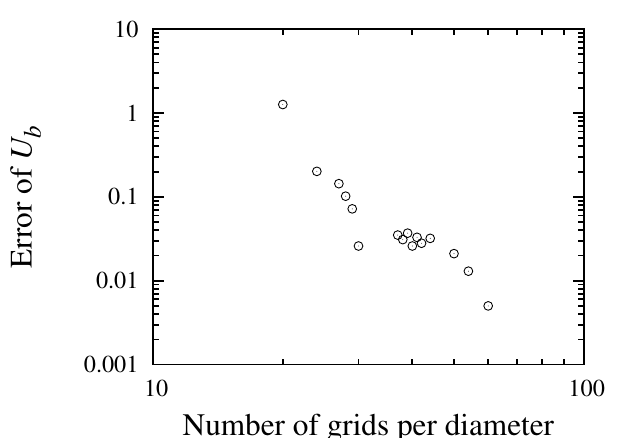}} 
 \end{center}
 \caption{(a) Mean pore velocity through a dense sphere pack as a function of number of grid cells per sphere diameter $D$. (b) Error 
    of mean pore velocity as a function of grid cells per diameter $D$.}
\end{figure}

We investigated the number of cells needed per sphere diameter for the bulk velocity to converge. Figure \ref{fig:pore-velocity-a} shows the intrinsic velocity versus the number of grid cells per diameter of the grains and Figure \ref{fig:pore-velocity-b} shows the logarithm of the error ($\varepsilon$) in the computed intrinsic velocity as a function of the logarithm of grid cells per diameter, taking the intrinsic velocity calculated using 70 grid cells per diameter as reference,

\begin{equation}  
\varepsilon = \frac{U_{i}-U_{i,ref}}{U_{i,ref}}.
\end{equation}

The  intrinsic velocity  converges monotonically with more than 35 cells per diameter and the error is limited to less than $5 \%$. There is no constant convergence rate due to the IBM method. However, on average, the convergence rate is at least of second order (Figure  \ref{fig:pore-velocity-b}).

\begin{figure}
 \begin{center}
  \includegraphics[width=.475\textwidth]{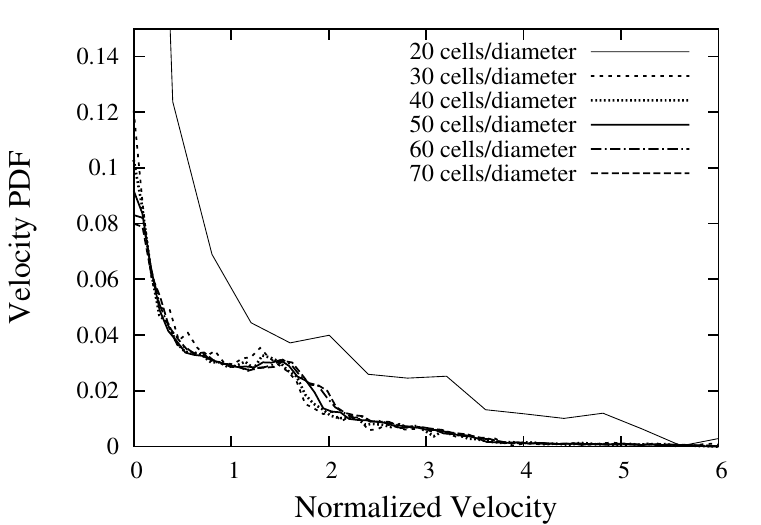}
  \caption{Velocity PDFs in a dense sphere pack for various number of grid cells per sphere diameter $D$.}
  \label{fig:pdf-resolution-study}
 \end{center}
\end{figure}

We inspect as well the probability density function (PDF) of the local velocity in pore space at different grid resolutions (Figure \ref{fig:pdf-resolution-study}). With more than 30 grid cells per diameter, the PDFs show only little variation.  We concluded that with 40 grid cells per sphere diameter it will be possible to get a sufficiently accurate velocity field and chose such a mesh for the simulations presented in this work. 

\paragraph{Flow through a random sphere pack.} 
We did further validation by comparing our results to the empirical correlations on Carman-Kozeny and Blake-Kozeny, respectively. Those relations make use of dimensional analysis to determine the overall form of the dependence of the permeability from porosity and grain diameter in a sphere pack, equation (\ref{eq:blake-kozeny}). The factor $\alpha$ in this relation is related to the ratio of the mean length of the passages a flow has to go through and the thickness of the layer that it goes through and is fitted to experimental measurements. Carman-Kozeny is connected to $\alpha=180$ while Blake-Kozeny is connected to $\alpha=150$,

\begin{equation}
\label{eq:blake-kozeny}
K=\frac{D^2 \epsilon^3}{\alpha(1-\epsilon)^2}\,.
\end{equation}

\begin{figure}
 \begin{center}
   \subfigure[]{\label{fig:sphere-pack}\includegraphics[width=.475\textwidth]{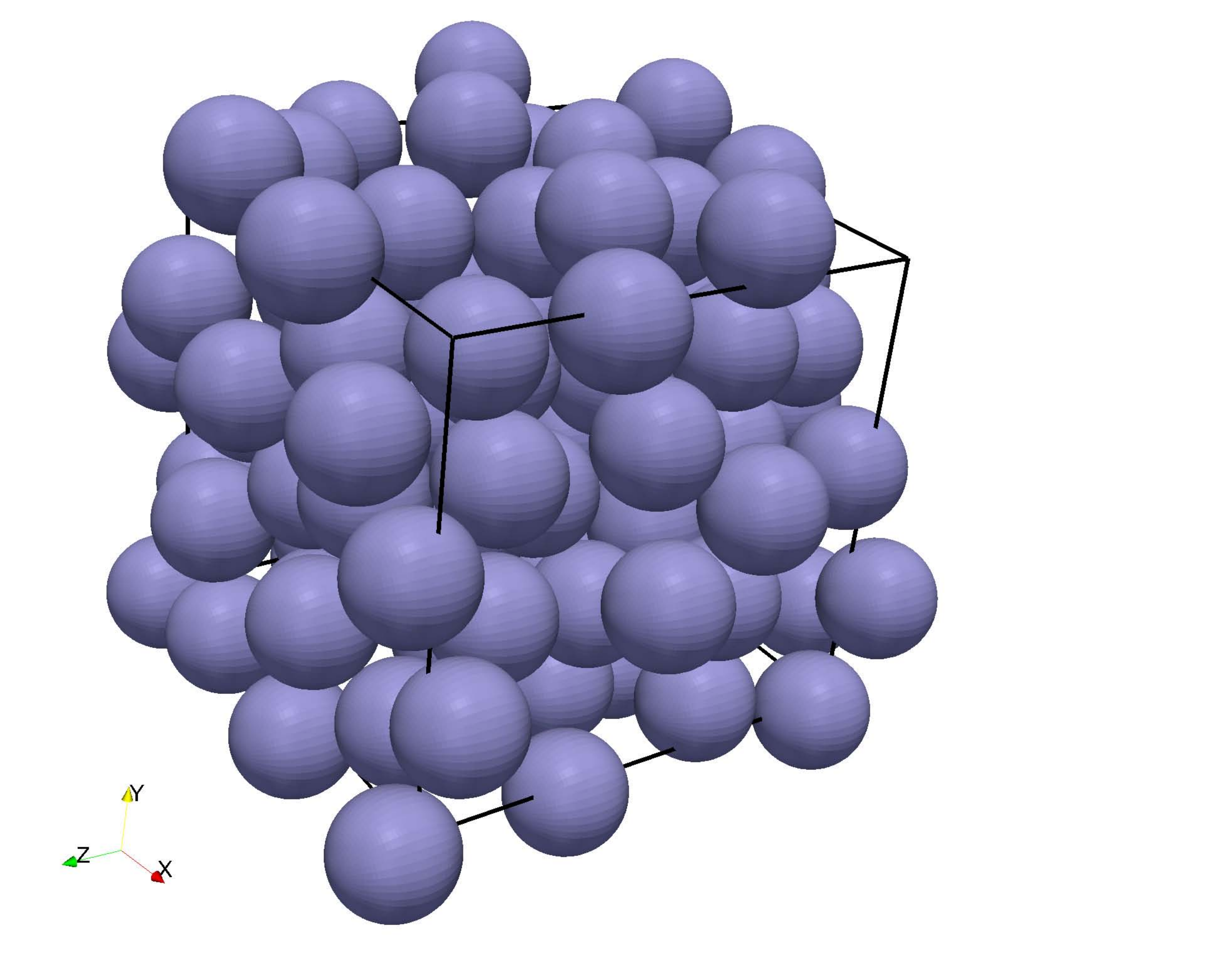}}
   \subfigure[]{\label{fig:blake-kozeny}\includegraphics[width=.5\textwidth]{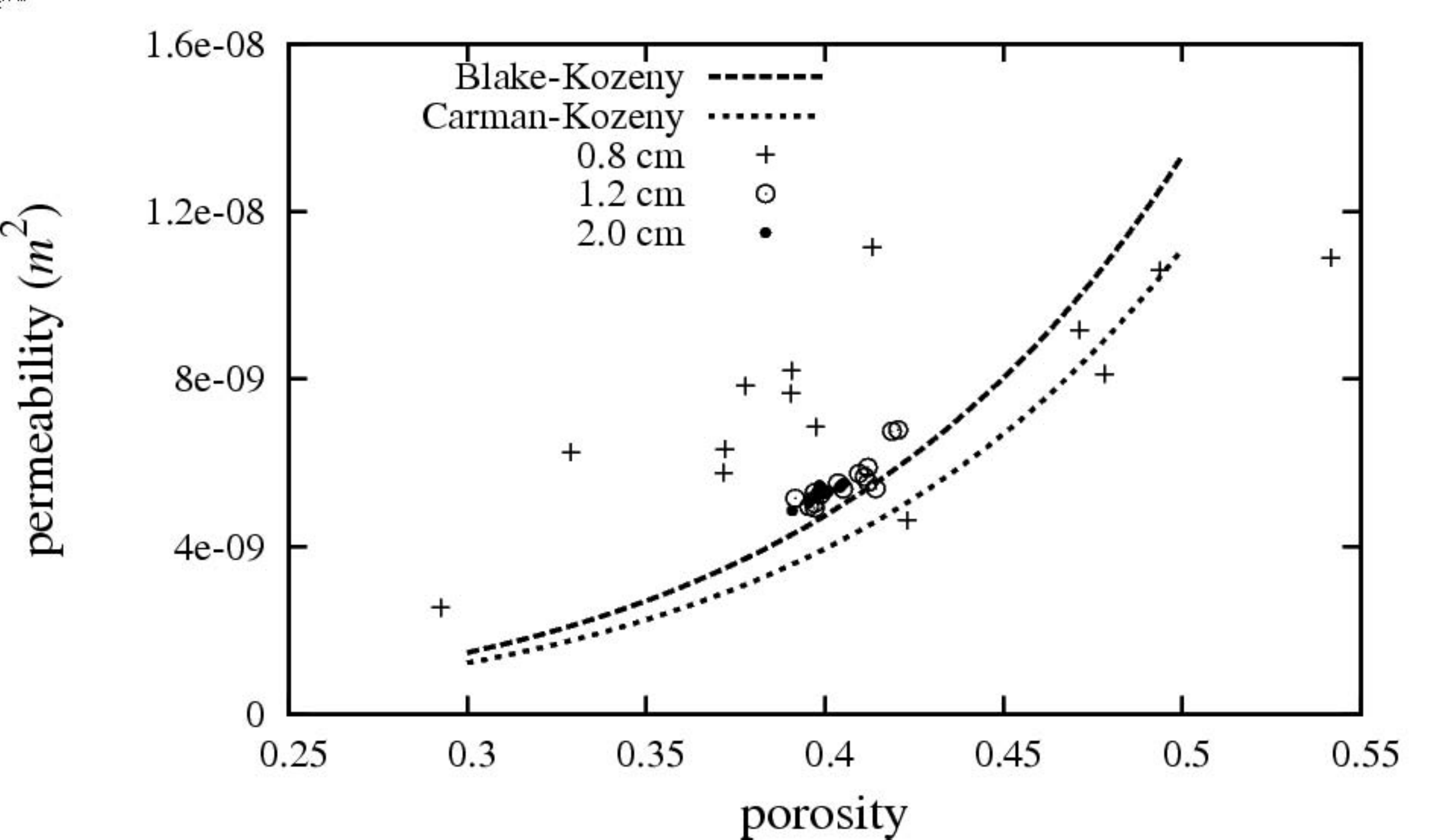}} 
 \end{center}
 \caption{(a) Random sphere pack: one realisation; 
          (b) Comparison of computed permeabilites in the inner domain of random sphere pack domains of  different sizes with the Blake-Kozeny and Carman-Kozeny relations.}
\end{figure}

A series of simulations through a random sphere pack with periodic boundary conditions in all three directions was conducted to find the minimum size of the REV. The grid resolution was $40$ cells per diameter. The size of the domain increased from $0.8\,\mathrm{cm}=4D$ to $2\,\mathrm{cm}=10D$. For each domain size, we simulated 15 different realisations of random sphere distributions, such as displayed in Figure \ref{fig:sphere-pack}, to obtain a reasonable sample size. By this series, we can check which domain size can be regarded as REV. We found that close to the domain boundaries our porosity was little larger than in the inner domain where it was distributed homogeneously. Therefore, we take only the values from the inner domain for comparison with the Blake-Kozeny relation. This inspection revealed that a domain size of $10D$ was sufficient to obtain in the inner region permeability values fully consistent with Blake-Kozeny's relation, see Figure \ref{fig:blake-kozeny}.

\begin{figure}
 \begin{center}
  \subfigure[]{\label{fig:PDF-inner-a}\includegraphics[width=0.475\textwidth] {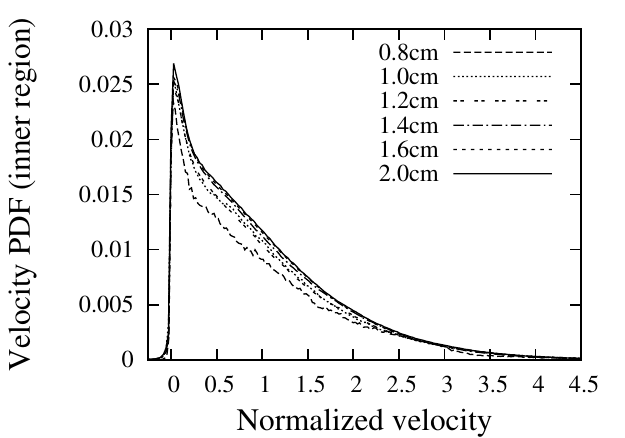}}
  \subfigure[]{\label{fig:PDF-outer-b}\includegraphics[width=0.475\textwidth] {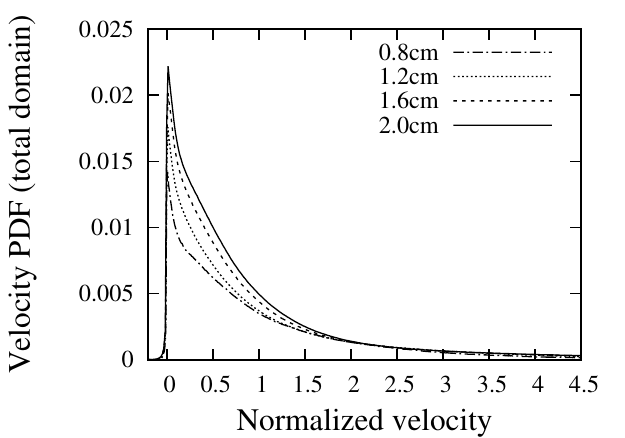}}
 
 \end{center}
 \caption{Velocity PDF in a random sphere pack at various domain sizes. (a) inner domain only; (b) total domain.}
\end{figure}

For every domain size we calculated the probability distribution function (PDF) of velocities in 
the range of $-2.6\times10^{-7}  \,\mathrm{m/s}$ and $8\times10^{-7}  \,\mathrm{m/s}$ using 1325 bins of size $8\times10^{-10}  \,\mathrm{m/s}$ for each realisation, and then averaged the PDF over all 15 realisations. Because of the varying porosity in our domains due to the special sphere packing procedure mentioned before, we first calculated the PDF of velocities of points residing in the 'inner region' of the domain only. That is to say, in each direction we omitted the points closer than 1.5 sphere diameters to the edge and then proceeded to calculate the PDF of the velocities as mentioned before. In the next step we calculated these PDFs for the complete domain too. The PDFs of the stream-wise velocity in random sphere packs is plotted in Figure \ref{fig:PDF-inner-a} for the inner domain and \ref{fig:PDF-outer-b} for the total domain. These plots demonstrate the convergence of the PDF with domain size. Here, the curves for the inner domain converge faster than the ones for the total domain which can be explained by the inhomogeneous porosity distribution close to the boundaries. The velocity in these plots is normalised by the intrinsic velocity $U_i$ which is the average velocity in the pore space. The distribution is highly skew. Maximum velocities of four times the averaged one can be observed, however with very small likelihood. 

Comparing our PDFs with those measured by magnetic resonance imaging \cite{deurer_04}, we observe large differences. Those were measured on various sample volumes, the smallest being in the range of sphere diameter. They represent velocities filtered on that scale. The maxima are at the order of magnitude of the pore velocity. Our PDFs have been evaluated at a sample size comparable to the grid spacing of the simulation which is much smaller than the sphere diameter. They can be regarded as unfiltered velocities and their maximum probability lies at values much smaller than the average pore velocity.

A striking feature of the PDFs are the negative velocities. Such negative velocities would not be expected in PDFs of the superficial velocities. They can be explained by the irregularity of the random sphere pack. This irregularity forces stagnation points at the front and back faces of the sphere to be off-centre. As a consequence, streamlines that travel to and from the stagnation points along the surfaces of the spheres have to point in negative x-direction in some regions and therefore generate negative stream-wise velocities. We conclude that those negative velocities can not be associated with flow separation in the traditional sense. Furthermore, we conclude that these negative velocities are not able to transport mass upstream over a long distance. However, they might increase the time a tracer needs to travel downstream and thus contribute to long tails of break-through curves.


\section{Parameter identification of an inhomogeneous permeability field.} 
In this section we focus on modeling flows in porous media on the macro-scale by the
Darcy equation (\ref{eq:darcy}). One key problem is to determine the averaged material properties,
here in particular the permeability tensors of the considered medium. Our approach 
is to determine them based on reference flow measurements taken from either
experiments or direct numerical simulation resolving the micro-scale behaviour of 
the media. Then the permeability tensors are chosen such that the resulting flow given 
by the Darcy equation for the experiment configuration matches the measurements 
optimally in a least-squares sense. Previous work on parameter estimation in 
similar settings includes~\cite{SchulzWittum1998}, \cite{MahnkenSteinmann2001} and~\cite{VexlerThesis2004}.

We outline an adjoint-based optimisation 
algorithm that performs the parameter fit for a suitable discretisation of the Darcy model. 
Special emphasis is on a discretisation for the problem which on the one hand satisfies
the necessary stability properties and on the other hand works well in the optimisation
context. Tests on some model configurations show the viability of the proposed method.

Our model for describing a fluid moving through a porous domain $\Omega\subseteq \mathbb R^d$
consists of the Darcy equation  (\ref{eq:darcy}) together with a volume integrated version of the mass balance equation~\eqref{eq:continuity}. After rearranging the Darcy equation, it reads
\begin{subequations}
  \label{eq:DarcyEq}
  \begin{align} 
    \Keffinv \uavg + \nabla \pavg &= 0, \\
    \nabla \cdot \uavg             &= f_p.
  \end{align}
\end{subequations}
The right hand side term $f_p$ is used to model sources and sinks within the
domain. By the position-dependent permeability tensor $\Keff\colon \Omega \to \mathbb R^{d\times d}$ we describe the 
effective permeability of the media at any given point in the domain. We use a tensor instead 
of a scalar quantity since not only isotropic but also anisotropic materials should be modelled.
According to \cite{Liakopoulos1965}, the tensor $\Keff$ is symmetric positive definite 
at any given point in $\Omega$.
For our test configurations we assume homogeneous Neumann boundary conditions and the condition 
$\int_\Omega f_p\, \text{d} x=0$ which ensures existence and uniqueness of solutions for suitably 
chosen spaces for velocity, pressure and the permeability tensor.

Due to the saddle point structure of~(\ref{eq:DarcyEq}), a finite element approximation has to 
be inf-sup stable. Since in the optimal control context we have to deal not only with the finite 
element spaces for the state variables but also with the corresponding dual spaces, using different Ansatz spaces
for pressure and velocity would add considerable complexity. Therefore
we want to use the same discrete spaces for both, pressure and velocity.
Hence the inf-sup condition has to be enforced by a suitable stabilisation. 
Here we use the local projection stabilisation (LPS)
approach (see~\cite{BraackSchieweck2011}) since compared to most other methods the 
resulting stabilisation terms are symmetric. Therefore the two approaches 
 \emph{``discretise-then-optimise''} and \emph{``optimise-then-discretise''} lead to the same 
set of discrete equations. In addition the systematic a posteriori error estimation approach
developed in~\cite{BeckerVexler2004} can be applied immediately. 
 For a detailed discussion
of LPS stabilisation for optimal control, see~\cite{Braack2009}. 
A LPS stabilised discretisation of the Darcy-Brinkman 
has been analysed in~\cite{BraackSchieweck2011}.
Their results include the Darcy equation with homogeneous isotropic media as a special case 
and can be extended in a straight-forward fashion towards non-homogeneous anisotropic media. 
We use bi-linear rectangular finite elements on a conforming grid which possesses a patch structure,
that is, the grid can be obtained by uniform refinement of a coarser grid $\mathcal{M}_h$.
Then the stabilised discretisation of~(\ref{eq:DarcyEq}) reads in weak form: find the discrete 
velocity and pressure $(\uavg_h,\pavg_h)$ which satisfy
\begin{align*}
&\qquad \int_\Omega \left \{ \Keffinv \uavg_h \varphi_v 
  - \pavg_h \nabla \cdot \varphi_v +\varphi_p \nabla \cdot \uavg_h \right \} \text{d}x \\
&\quad + \sum_{M\in\mathcal M_h} \int_{M} \left \{ h_M^2 \kappa_M(\nabla \uavg_h) \kappa_M(\nabla \varphi_v) 
+ \kappa_M (\nabla \pavg_h) \kappa_M(\nabla \varphi_p)
\right \}\text{d}x \\
&= \int_\Omega f_p \varphi_p \,\text{d}x
\end{align*}
for all discrete test functions $(\varphi_v,\varphi_p)$. 
The fluctuation operator $\kappa_M$ is defined locally on each cell $M\in\mathcal{M}_h$ 
of the coarser 
grid as $\kappa_M = \Id-\Pi_M$ with $\Id$ denoting the identity and $\Pi_M$ the $L^2$ 
projection onto the space of constant functions on $M$. The diameter of $M$ is denoted by $h_M$.
 Stability and first order convergence in the $L^2$ norm
 with respect to the discretisation parameter $h$ are shown in~\cite{Himmelstoss2011}.

For the parameter estimation problem we assume that we have a priori information about the distribution
of different materials within the domain, furthermore that 
the domain can be divided into finitely many sharply bounded 
regions with different materials and that within each region the effective permeability tensor stays 
constant. In order to avoid enforcing the positive definiteness of the permeability tensor 
by additional constraints, we parametrise $\Keffinv$ in a suitable way by a finite number of parameters 
$q_i \in \mathbb R$. If we restrict our considerations to materials with diagonal permeability tensor, 
then a possible parametrisation consists of the $d$ diagonal entries  of $\Keffinv$ on each
region. To ensure positive definiteness, the vector of parameters $q$ is bounded away from zero by algebraic 
constraints.

Computing $\uavg$ and $\pavg$ given $\Keffinv(q)$ is a well-posed problem, however the inverse problem of
determining $q$ from given measurements of $\uavg$ and $\pavg$ can be ill-posed, that is, small variations 
in the measurement data can lead to big variations in the recovered $q$. Therefore  we apply a Tychonoff
regularisation with parameter $\alpha\geq 0$ such that the parameter identification problem can be 
stated as: Minimise 
\begin{align}
  \label{eq:CostFunctional}
  J(q,u) &= \|C u - z \|^2  + \frac\alpha2 |q |^2
\end{align}
subject to $u=(\uavg,\pavg)$ solving the Darcy equation~(\ref{eq:DarcyEq}) for $\Keffinv := \Keffinv(q)$
and $q\in Q_{\text{ad}} \subseteq Q=\mathbb R^{N}$ where $N$ is the number of parameters in the parametrisation 
of the permeability tensor. The linear operator $C$ models some measurements done on the computed solution,
this could be for example evaluation of the velocity field at certain points within the domain. The 
value $z$ represents the corresponding reference data for that measurement obtained from a micro-scale
model or from an experiment. Since from a micro-scale simulation in principle we can obtain a 
complete reference state, it makes sense to chose the identity as observation operator $C$ in that case. 
The parameter identification problem can be interpreted as an optimal control
problem where the control variable $q$ should be chosen in such a way that the state variable $u$ matches a 
desired state described by the measurements as good as possible. We enforce positive definiteness of the 
permeability tensor by an appropriate choice of the closed set $Q_{\text{ad}}\subseteq Q$.

The existence of a solution to the optimal control problem can be shown by standard arguments, see for example
the textbook~\cite{Troeltzsch2010}. Since the problem is in general non-convex, uniqueness of the solution
cannot be guaranteed without further assumptions.

As noted before, for any control $q$ there is a unique state $u$ satisfying~(\ref{eq:DarcyEq}). Therefore 
we can define the control-to-state mapping 
\begin{equation*}
  S\colon q \mapsto u
\end{equation*}
with $u=(\uavg,\pavg)$ solving the Darcy equation~(\ref{eq:DarcyEq}) for $\Keffinv := \Keffinv(q)$. We introduce the reduced 
cost functional $j(q) := J(q,S(q))$ and state the reduced optimisation problem
\begin{equation*}
  \min   j(q) \text{ subject to $q\in Q_{\text{ad}}$.}
\end{equation*}
To solve this reduced problem we use a primal-dual-active-set strategy (PDAS) 
(see, e.\,g., \cite{HintermItoKunisch2003}) 
to treat the algebraic constraints on 
$q$ resulting from the choice of $Q_{\text{ad}}$. In each step of the PDAS, an unconstrained optimisation problem 
has to be solved. For that purpose a globalised Newton-CG method is used. Gradient and Hessian information
are computed via an adjoint approach, for further details on the algorithm see, 
e.\,g.~\cite{BeckerMeidnerVexler2007} or~\cite{VexlerThesis2004}. To ensure fast convergence of the
Newton method, exact derivatives that are consistent with the discrete stabilised state equation are 
essential. Therefore in particular the derivatives of the stabilisation terms with respect to $q$ have
to be taken into account when deriving the auxiliary equations used for Hessian evaluation. 

Considering the computational complexity of the outlined algorithm, we note that the number of Newton steps does not depend on the fineness of the discretisation. The inner CG solver takes in the worst case $\mathcal O(N)$ iterations and for each
iteration we have to solve two auxiliary PDEs, which each take $\mathcal O(L)$ operations with a multi-grid 
solver, where $L$ is the dimension of the finite element space. So in total we expect our algorithm to have the
complexity $\mathcal O(N\cdot L)$.

\begin{figure}
  \centering
  \def\svgwidth{0.9\textwidth}
  \includegraphics[width=\svgwidth]{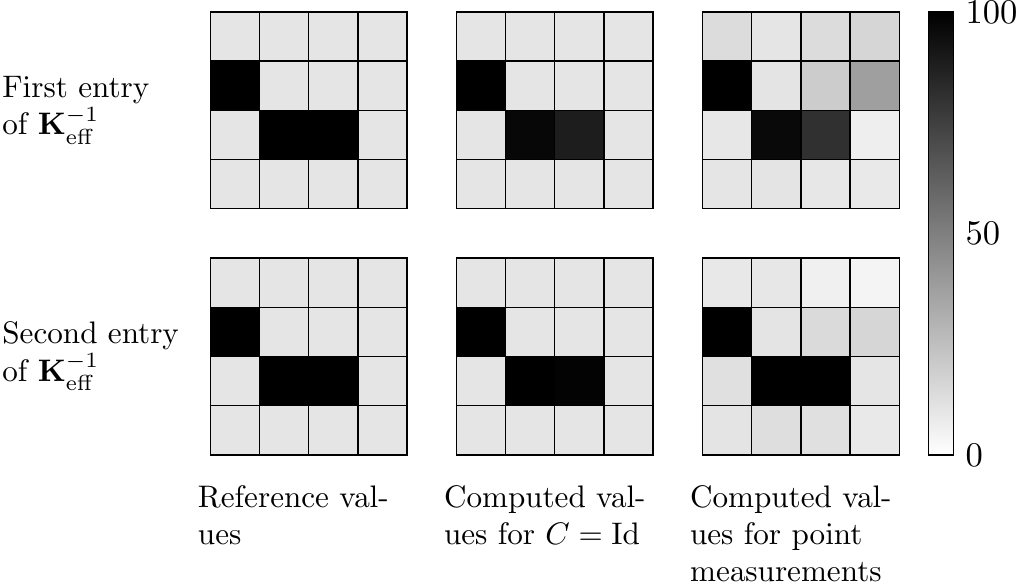}
  \caption{Values of $\Keffinv$ over the domain $\Omega$ for $C=\Id$ and for 32 point measurements}
  \label{fig:ParameterEstTensors}
\end{figure}

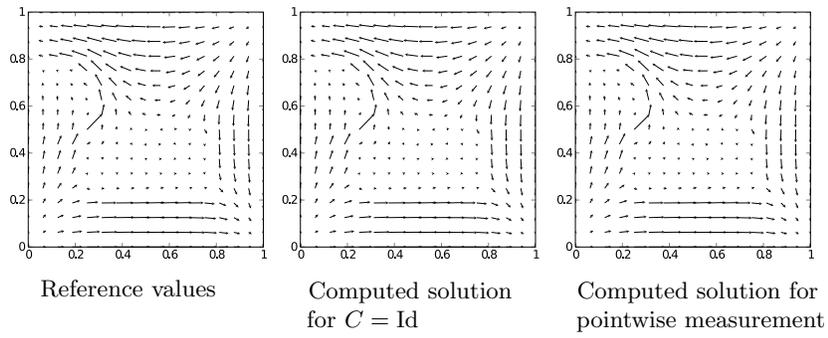
\begin{figure}
  \centering
  \def\svgwidth{1.\textwidth}

\begingroup
  \makeatletter
  \providecommand\color[2][]{
    \errmessage{(Inkscape) Color is used for the text in Inkscape, but the package 'color.sty' is not loaded}
    \renewcommand\color[2][]{}
  }
  \providecommand\transparent[1]{
    \errmessage{(Inkscape) Transparency is used (non-zero) for the text in Inkscape, but the package 'transparent.sty' is not loaded}
    \renewcommand\transparent[1]{}
  }
  \providecommand\rotatebox[2]{#2}
  \ifx\svgwidth\undefined
    \setlength{\unitlength}{450.85712891bp}
    \ifx\svgscale\undefined
      \relax
    \else
      \setlength{\unitlength}{\unitlength * \real{\svgscale}}
    \fi
  \else
    \setlength{\unitlength}{\svgwidth}
  \fi
  \global\let\svgwidth\undefined
  \global\let\svgscale\undefined
  \makeatother
  \begin{picture}(1,0.44740178)
    \put(0,0){\includegraphics[width=\unitlength]{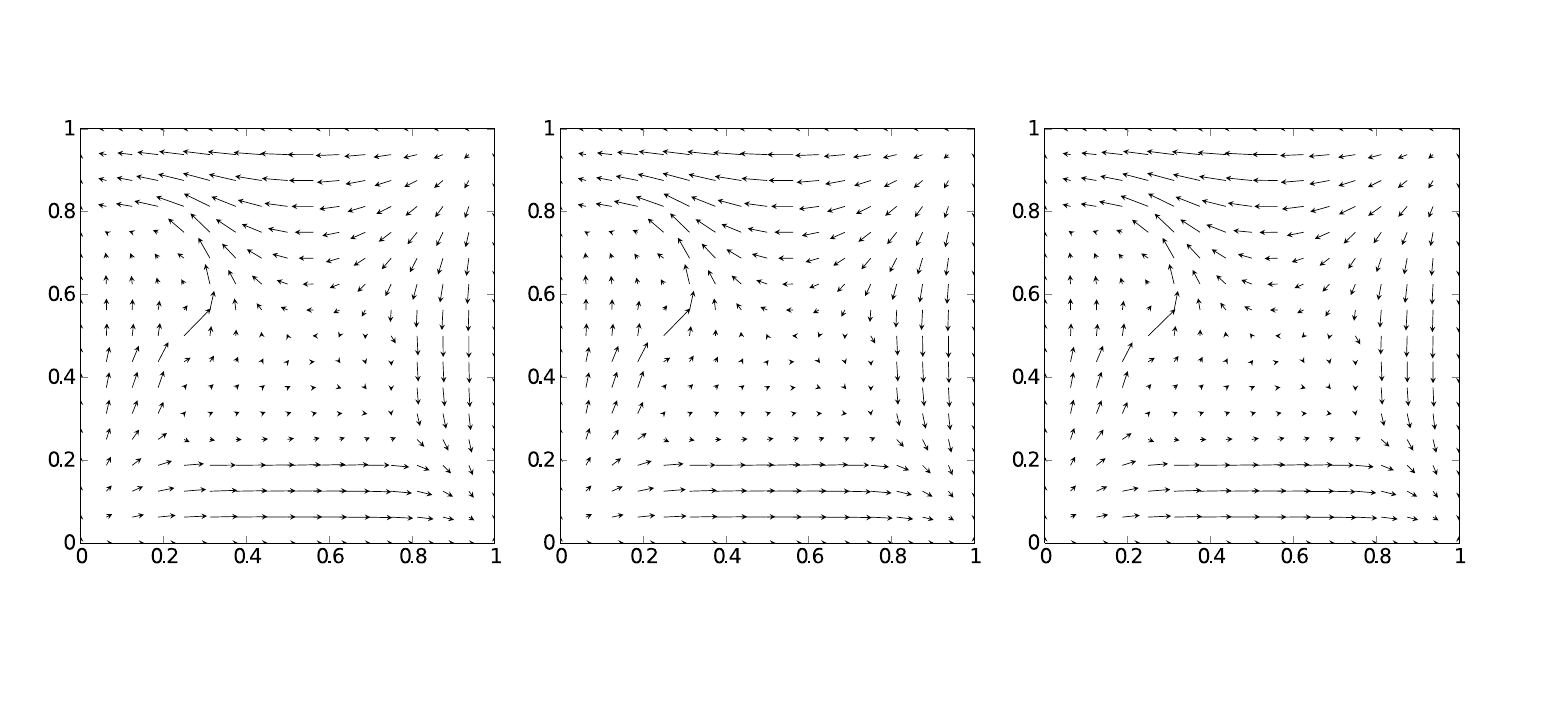}}
    \put(0.0650395,0.04489973){\color[rgb]{0,0,0}\makebox(0,0)[lb]{\smash{Reference values}}}
    \put(0.66934879,0.04236487){\color[rgb]{0,0,0}\makebox(0,0)[lb]{\smash{Computed solution for}}}
    \put(0.36690558,0.04236482){\color[rgb]{0,0,0}\makebox(0,0)[lb]{\smash{Computed solution}}}
    \put(0.66869455,0.01064637){\color[rgb]{0,0,0}\makebox(0,0)[lb]{\smash{pointwise measurement}}}
    \put(0.36527251,0.01064637){\color[rgb]{0,0,0}\makebox(0,0)[lb]{\smash{for $C=\Id$}}}
  \end{picture}
\endgroup

  \caption{Comparison of exact velocity field and velocity fields resulting from estimated $q$}
  \label{fig:VelocityFields}
\end{figure}

For the numerical tests we consider the Darcy problem on the two-dimensional unit square $\Omega=(0,1)^2$.
We subdivide $\Omega$ into 16 equally sized squares $\Omega_i$, $i=1,\dots,16$ 
and assume that on each square the permeability tensor 
is constant and can be represented by a diagonal matrix. Therefore we choose the control space 
$Q=\mathbb R^{32}$ and define the parametrisation of the permeability tensor by
\begin{equation*}
  \Keffinv (q) |_{\Omega_i} = 
  \left (  
    \begin{array}{cc}
      q_{2i-1} & 0 \\
      0 & q_{2i}  \\
    \end{array}
  \right ) \quad \text{for $i=1,\dots,16$.}
\end{equation*}
For convenience we denote the vector collecting all the entries in the first component of $\Keffinv$
by $q^A \in \mathbb R^{16}$ and the one collecting the entries in the second component by $q^B$
The source term is chosen as 
\begin{equation*}
  f_p(x,y) = 2\cos(\pi x)\cos(\pi y),
\end{equation*}
and the set of admissible controls is defined as
\begin{equation*}
  Q_{\text{ad}} = \left \{
    q\in \mathbb R^{32}
    \middle |
     q \geq 1
  \right \}.  
\end{equation*}
Since the problem is reasonably well conditioned, we can omit the regularisation term by setting $\alpha=0$.  
For the discretisation of pressure and velocity, a grid with 4096 cells is used. The measurement data 
$z$ are generated synthetically
by performing a forward simulation with a reference parameter vector $q_{\text{ref}}$. We investigate two
choices for the observation operator $C$, first the identity and second an operator modelling 
32 point measurements of pressure and velocity within the domain. 
A visual comparison of the reference permeability tensor and the permeability tensors computed by the
parameter identification algorithm can be seen in Figure~\ref{fig:ParameterEstTensors}. For both choices
of the observation operator $C$, good qualitative agreement between the reference and the 
computed permeability values is observed. However, for the case  $C=\Id$, the estimated parameters
 are better than for the point-wise measurements 
since more data enters the computation. These observations are confirmed when looking at the relative errors
$\frac{ \|q^A - q^A_{\text{ref}} \|_2 }{\| q^A_{\text{ref}} \|_2 }$ and $\frac{ \|q^B - q^B_{\text{ref}} \|_2 }{\| q^B_{\text{ref}} \|_2 }$
listed in Table~\ref{tab:ParameterErrors}.
A qualitative comparison of the resulting velocity fields to the 
reference velocity field is shown in Figure~\ref{fig:VelocityFields}.

\begin{table}
  \centering
  \caption{Relative errors of the two tensor components for both choices of $C$}
  \begin{tabular}{lrr}
     &\qquad $C=\Id$ & $\qquad $  point measurements \\
    \hline
    $\frac{ \|q^A - q^A_{\text{ref}} \|_2 }{\| q^A_{\text{ref}} \|_2 }$ & 0.0655 & 0.181 \\
    $\frac{ \|q^B - q^B_{\text{ref}} \|_2 }{\| q^B_{\text{ref}} \|_2 }$ & 0.00565 & 0.0634 \\
  \end{tabular}
  \label{tab:ParameterErrors}
\end{table}

\clearpage


\section{High order Finite Element Method for the advection diffusion equation.} 
One of the main numerical problems for simulations of tracer transport on the macro-scale are strong gradients within the tracer fields developing in situations were convection dominates over diffusion. 
Standard Bubnov-Galerkin finite elements are known to deliver oscillating solutions for convection dominated problems for meshes which are not fine enough. It has not yet been proved whether raising the polynomial degree of the shape functions will increase or decrease numerical oscillations. This paper will show that an increase of the polynomial degree ($p$-FEM) stabilises the numerical oscillations in Bubnov-Galerkin type finite elements naturally without adding any additional stabilisation term.

We will demonstrate the improvement of the numerical accuracy with polynomial order using a one-dimensional stationary convection-diffusion problem (\ref{eq:convection-diffusion-averaged-2}). Given a constant convection velocity $u_x$, a steady and constant effective diffusion coefficient $\Gamma^{\text{eff}}$ and a source term $f$, the problem is to find $c:\Omega\rightarrow\mathbb{R}$, such that with Dirichlet boundary conditions

\begin{equation}
\label{eq:pfem_cd_strong_1d}
\left\{
\begin{aligned}
u_x \frac{dc}{dx} - \Gamma^{\text{eff}} \frac{d^2 c}{dx^2} &= f \quad &\text{on}& \quad \Omega = \{x | 0 < x < 1\} \\
c &= 0 &\text{at}& \quad x=0 \\
c &= 0 &\text{at}& \quad x=1 \\
\end{aligned}
\right.
\end{equation}

We contrast the numerical errors of $p$-FEM \cite{Szabo:04.1} to the standard $h$-FEM \cite{Szabo:91} in which linear shape functions are used and follow the analysis scheme presented in \cite{Donea:03}. Herein, the truncation error of a Bubnov-Galerkin discretisation is quantified in order to specify the additional diffusion term used in Petrov-Galerkin methods. For a $h$-FEM, this results in a discretised equation which includes the numerical diffusion $\bar{\Gamma^{\text{eff}}}$
 
\begin{equation}
\label{eq19}
u_x \left(\dfrac{c_{j+1}-c_{j-1}}{2h}\right)-(\Gamma^{\text{eff}}+\bar{\Gamma^{\text{eff}}})\left(\dfrac{c_{j+1}-2c_j+c_{j-1}}{h^2}\right)=1\,.
\end{equation}

The extra term $\bar{\Gamma^{\text{eff}}}$ can be interpreted either as the truncation error of the Bubnov-Galerkin method of first order or as an additional diffusivity required to provide nodally exact results. This term is a function of the mesh $P\acute{e}clet$ number and reads

\begin{equation}
\label{eq20}
\bar{\Gamma^{\text{eff}}}=\left(\coth{Pe}-\dfrac{1}{Pe}\right)\Gamma^{\text{eff}} Pe\,.
\end{equation}

The mesh $P\acute{e}clet$ number is defined as

\begin{equation}
Pe=\dfrac{u_x h}{2\Gamma^{\text{eff}}}\,.
\label{fig:mpec}
\end{equation}

where $h$ is the mesh or grid size.

The value of $\bar{\Gamma^{\text{eff}}}$ increases with the mesh $P\acute{e}clet$ number. In fact, equation~(\ref{eq20}) forms the basic motivation behind using the Petrov-Galerkin method. In many stabilisation approaches, one tries to control the artificial numerical oscillations in convection dominated problems by compensating for the truncation error by means of adding additional diffusivity. However, it will be shown in the next section that the truncation error of the Bubnov-Galerkin method is decreased by a mere increase of the polynomial order of the spatial discretisation.

It is important to mention here that the truncation error study shown in next sections is also performed in more details in \cite{cai_13}. In the paper \cite{cai_13}, the stabilization capability of the \textit{p}-FEM for convection-dominated transport problems is explained mathematically by analyzing stiffness matrices. Numerical examples show that using sufficiently high order polynomial degrees for shape functions can eliminate the nodal oscillations in numerical solutions for convection-dominated problems, where the mesh $P\acute{e}clet$ number is greater than one. This approach will be introduced in following sections again in order to explain, why the high order FEM is suitable for solving convection-dominated problems of tracer transport on the macro-scale.

\subsection{Truncation error of the Bubnov-Galerkin discretisation in the $p$-FEM}
\label{sec:pFEM_error_pfem}

In this section, the truncation error of $p$-FEM for the same example as presented above is considered, where hierarchic shape functions of second order derived from the set of integrated Legendre polynomials are applied and the polynomial orders up to 5 are investigated. Compared to Lagrange shape functions, hierarchic shape functions are easy to construct since lower order shape functions are subsets of higher order ones. We refer to \cite{Szabo:91}, where the complete hierarchy of spaces is introduced. 

In general, the system equation using polynomial degrees higher than 2 can be also condensed analogously as in equation~(\ref{eq19}), using $\bar{\Gamma^{\text{eff}}}_p$ instead of $\bar{\Gamma^{\text{eff}}}$ as all higher modes are purely internal to the element. \\

Analogous to the previous analysis, one can get the following diffusion using second to fifth order polynomials for shape functions, respectively.
\begin{equation}
\label{eq30}
\begin{aligned}
\bar{\Gamma^{\text{eff}}}_2 &=\dfrac{1}{3}Pe^2\Gamma^{\text{eff}} \\
\bar{\Gamma^{\text{eff}}}_3 &=\dfrac{5Pe^2\Gamma^{\text{eff}}}{Pe^2+15} \\
\bar{\Gamma^{\text{eff}}}_4 &=\dfrac{\Gamma^{\text{eff}}(Pe^4+35Pe^2)}{10Pe^2+105} \\
\bar{\Gamma^{\text{eff}}}_5 &=\dfrac{14\Gamma^{\text{eff}}(4Pe^4+90Pe^2)}{4Pe^4+420Pe^2+3780}
\end{aligned}
\end{equation}

\begin{figure}[ht]
\centering
\includegraphics[angle=90,scale=0.5]{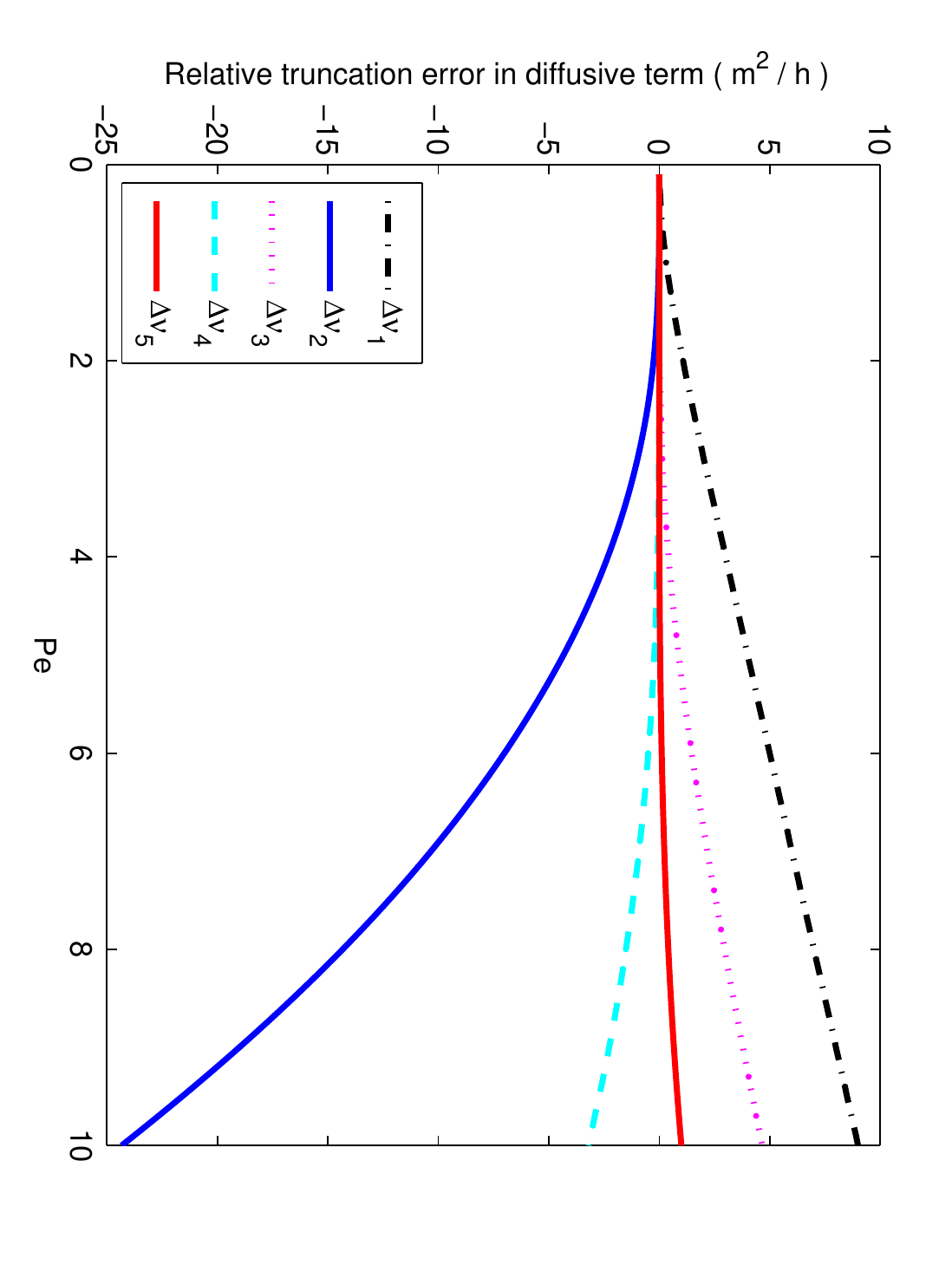}
\caption{Truncation error with different polynomial degrees}
\label{fig:DeltaNuP}
\end{figure}

The truncation error of \textit{p}-FEM is defined as

\begin{equation}
\label{eq29}
\Delta{\Gamma^{\text{eff}}_p}=\bar{\Gamma^{\text{eff}}}-\bar{\Gamma^{\text{eff}}}_p
\end{equation}

and depicted in dependence of $Pe$ in Figure~\ref{fig:DeltaNuP}, where the ordinate displays $\Delta\Gamma^{\text{eff}}_p$.\\

In general, the curves have different tendencies which correspond to the parity of the polynomial degree. Odd degrees generate curves which increase monotonically as $Pe$ increases while the even ones decrease. Although the sign of truncation error depends on the parity of the order, the absolute value of truncation error decreases when the order of shape functions grows. Accordingly, the numerical solution at nodes approaches the exact solution. \\
On the other hand, using odd polynomial degrees, the numerical diffusivity of the high order approach is less than $\bar{\Gamma^{\text{eff}}}$. This lack of diffusivity is the reason of the oscillatory behaviour of the numerical solution at high $Pe$. By contrast, the numerical diffusivity is always greater than $\bar{\Gamma^{\text{eff}}}$ using even polynomial degrees. Consequently, nodal solutions exhibit an over-diffusive behaviour and never show nodal oscillations. This result is further analysed from a mathematical perspective in the next section.\\

\subsection{Connection of the stability and the structure of the system matrix}
\label{sec:pFEM_stab_matrix}

Stability, i.e. oscillations or not, is determined by the structure of the system matrix. The numerical simulation will start to oscillate if the discrete maximum principle is violated \cite{Rank:83}. Considering a system matrix structure such as given in equation (\ref{eq:tridiag_alpha}), it can be proved that no oscillations occur for $\alpha<1$ \cite{Ernst:00}.  
\begin{equation}
\label{eq:tridiag_alpha}
\mathcal{A}(\alpha) = tridiag(-1-\alpha, 2, -1+\alpha)
\end{equation}

The system matrix resulting from the condensed equation (\ref{eq19}) can be written as 
\begin{equation}
\label{eq:alpha_p}
\begin{aligned}
\mathcal{A}_{p} = \dfrac{(\Gamma^{\text{eff}}+\bar{\Gamma^{\text{eff}}}_p)}{h^2}tridiag(-1-\alpha_p, 2, -1+\alpha_p)\,, \\
\alpha_p=\dfrac{u_x h}{2(\Gamma^{\text{eff}}+\bar{\Gamma^{\text{eff}}}_p)}\,.
\end{aligned}
\end{equation}

\begin{figure}[ht]
\centering
\includegraphics[width=0.9\textwidth]{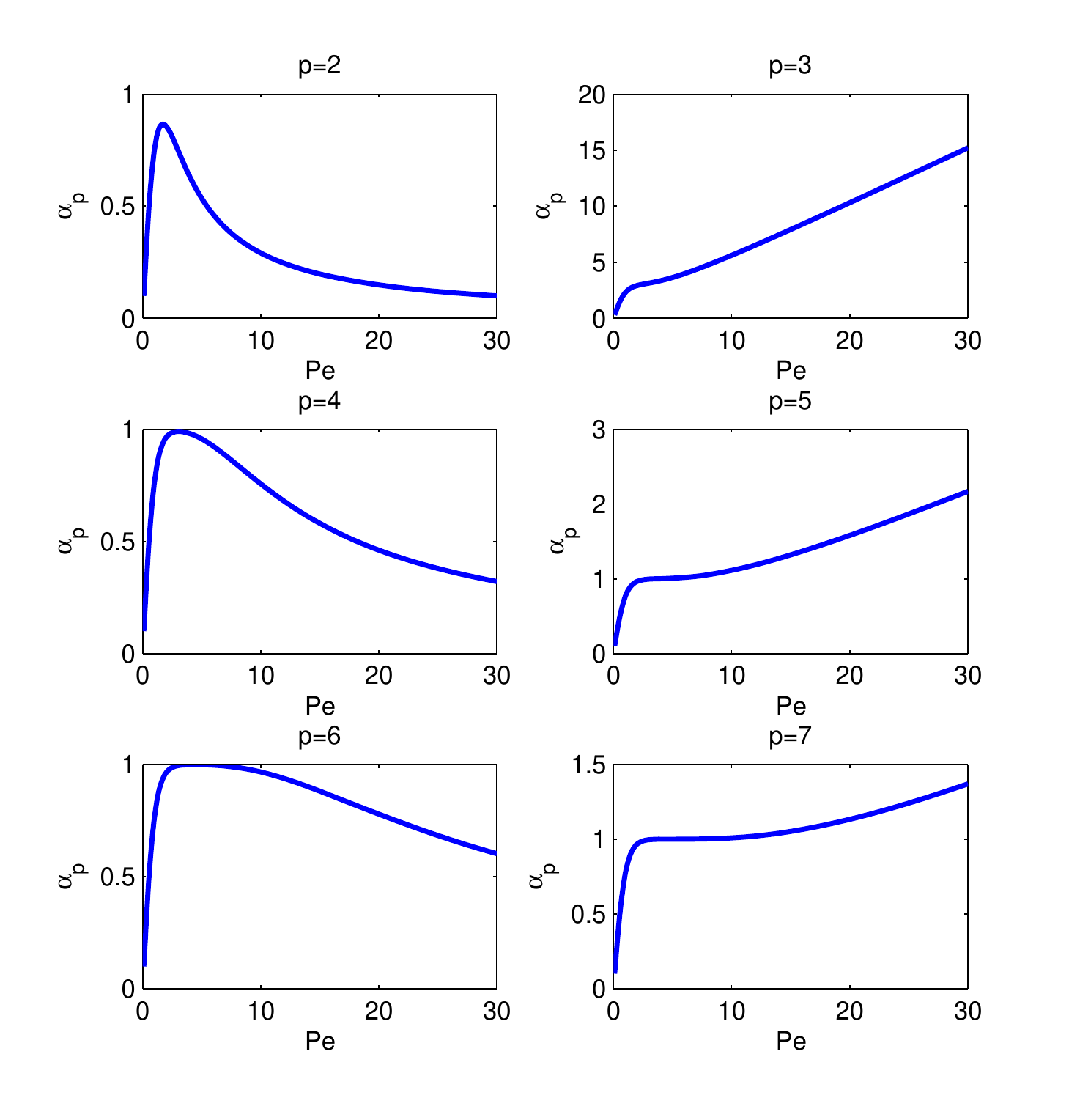}
\caption{$\alpha_p$ behaves differently for odd and even polynomial degrees}
\label{fig:alpha}
\end{figure}

\begin{figure}[ht]
\includegraphics[width=0.9\textwidth]{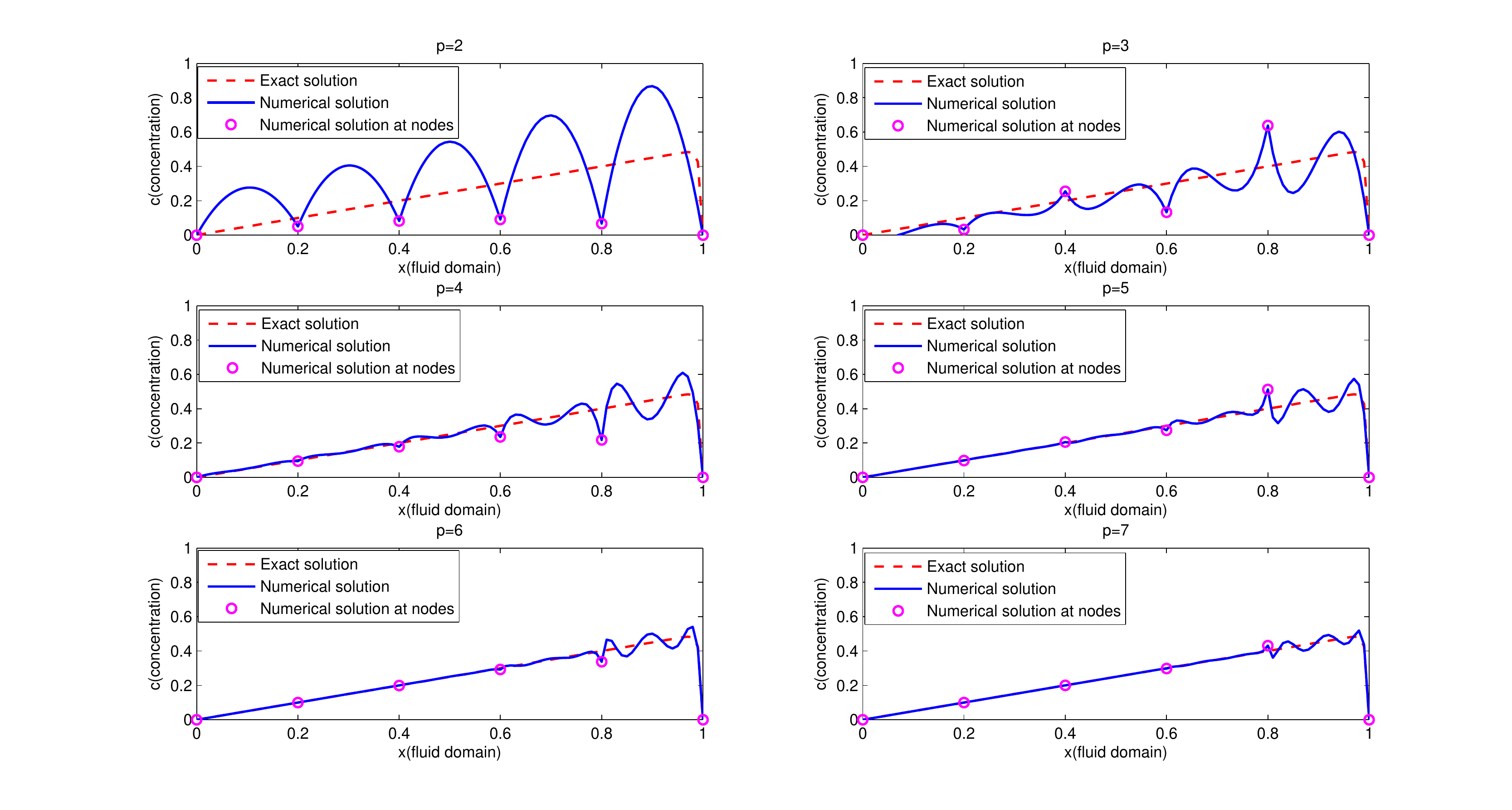}
\caption{Comparison of numerical, exact and nodal solutions with different Ansatz degree, Pe=20}
\label{fig:sol_p}
\end{figure}

Consequently, the stability of nodal solutions is determined by the value of $\alpha_p$. Further, the value of $\alpha_p$ can be quantified for higher order polynomial degrees based on equation~(\ref{eq30}):

The corresponding values are plotted in Figure~\ref{fig:alpha}. It can be observed that for odd polynomial degrees $\alpha_p$ increases as $Pe$. For even polynomial degrees, $\alpha_p$ first increases and then decreases while the value is always smaller than 1. This in turn means that for even polynomial degrees, the numerical solution at nodal degrees of freedom never oscillates. This result also coincides with the conclusion from the truncation error analysis in the previous section. To further clarify this point, we plot the solution of the 1D example with $Pe=20$ shown in Figure~\ref{fig:sol_p}.\\

 Here, the exact solution denotes the analytical solution of the differential equation~(\ref{eq:pfem_cd_strong_1d}). Figure~\ref{fig:sol_p} illustrates that when the polynomial degree is even, numerical oscillations only stem from internal modes and numerical solutions at each node do not oscillate. For odd polynomial degrees, numerical oscillations are reflected by both internal and nodal degrees of freedom.\\

\begin{figure}[ht]
\centering
\includegraphics[width=0.5\textwidth]{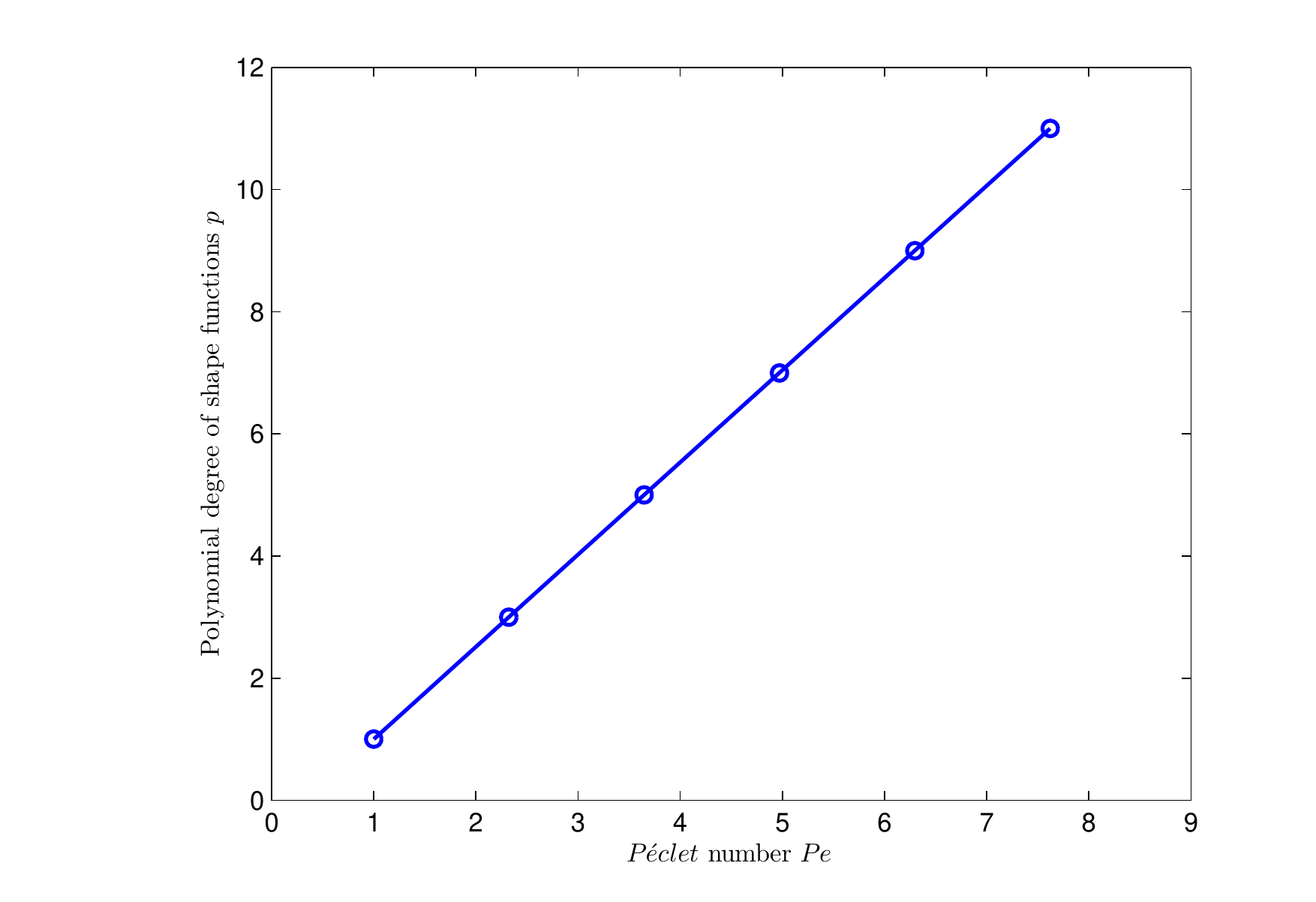}
\caption{The relation between a given mesh $P\acute{e}clet$ number and the minimum required polynomial degree}
\label{fig:pep}
\end{figure}

By setting $\alpha_p=1$ in equation~(\ref{eq:alpha_p}), we can compute the maximum allowed $Pe$ which guarantees nodally stable solutions for the given polynomial degree of the shape functions.
In other words, for a given mesh $P\acute{e}clet$ number, the corresponding $p$ stated in equation~(\ref{eq:oddPe}) is the minimum required polynomial degree and their relationship is depicted in Figure~\ref{fig:pep}. It turns out to be almost linear for polynomial orders $p\leq11$. 

\begin{equation}
\label{eq:oddPe}
\begin{aligned}
p=3 &  & Pe = 2.322185 \\
p=5 &  & Pe = 3.646738 \\
p=7 &  & Pe = 4.971786 \\
p=9 &  & Pe = 6.297019 \\
p=11 &  & Pe = 7.622340 \\
\cdots
\end{aligned}
\end{equation}

Different from other up-winding methods, where the additional efforts for modelling the necessary artificial diffusivity for more complicated problems, the high order FEM generates the additional numerical diffusion naturally by purely increasing the polynomial degrees. In the following example, numerical results of the one-dimensional convection-diffusion transport problem are compared to the exact solution. The given differential equation~(\ref{eq:pfem_eg_1d_single})

\begin{equation}
\label{eq:pfem_eg_1d_single}
\left\{
\begin{aligned}
a\frac{dc}{dx} - \Gamma^{\text{eff}} \frac{d^2 c}{dx^2} &= 0 \quad &\text{on}& \quad \Omega = \{x | 0 < x < 1\} \\
c &= 0 &\text{at}& \quad x=0 \\
c &= 1 &\text{at}& \quad x=1
\end{aligned}
\right.
\end{equation}

has the analytical solution

\begin{equation}
\label{eq:pfem_eg_1d_single_anal}
y = \frac{e^{ax/\Gamma^{\text{eff}}} - 1 }{ e^{a/\Gamma^{\text{eff}}} - 1 } \, .
\end{equation}

\begin{figure}[ht]
\begin{center}
\includegraphics[width=.9\textwidth]{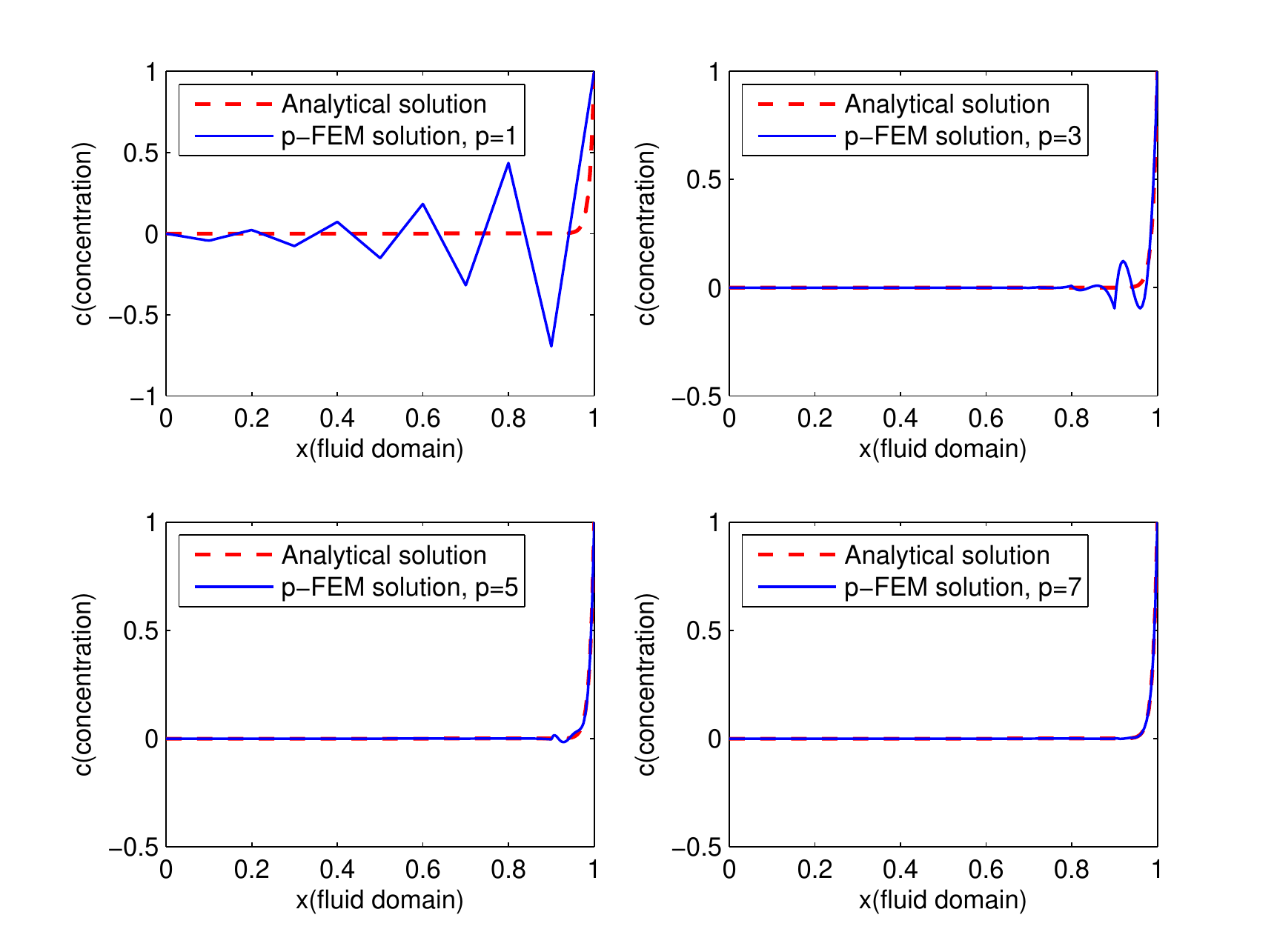}
\caption{Numerical solution with different polynomial degrees, $Pe=5$}
\label{fig:pfem_eg_1d_single}
\end{center}
\end{figure}

When the mesh is fixed, the ratio between a velocity and a diffusivity determines the mesh $P\acute{e}clet$ number and characterises the convergence of the numerical solution. When the mesh $P\acute{e}clet$ number increases, the standard Bubnov-Galerkin method based on linear elements exhibits oscillations in the numerical solution. We choose the parameters $a=2\,\mathrm{m/h}$, $\Gamma^{\text{eff}}=0.02\,\mathrm{m^2/h}$, and compute the corresponding numerical solutions with 10 elements of the same length $h=0.1$. Figure~\ref{fig:pfem_eg_1d_single} shows numerical solutions with different polynomial degrees. The dashed line denotes the exact solution while the solid line represents the numerical solution.\\

As expected, when mesh $P\acute{e}clet$ number $Pe=\dfrac{ah}{2\Gamma^{\text{eff}}}=5$ is larger than 1, the numerical solution with linear Bubnov-Galerkin discretisation introduces non-physical oscillations. The $p$-FEM can eliminate these oscillations by simply raising the polynomial degree $p$. It is observed in Figure~\ref{fig:pfem_eg_1d_single} that with $p=7$, the oscillation is drastically suppressed and the numerical solution is in good agreement with the analytical one.



\section{Conclusions}
In this paper we presented some efforts to improve understanding and simulation of flow and transport in porous media. Using consistent volume averaging, it can be shown that traditional closures, such as effective permeability and diffusivity are not applicable in all situations. Those situations arise for dispersion in the initial phase of tracer transport, for strongly inhomogeneous permeability fields and for convection dominated transport.

The initial phase of tracer transport is characterised by non-Gaussian tracer plumes, the so-called non-Fickian regime. The transition from non-Fickian to Fickian dispersion is dependent on how long tracer patches stay in low-speed regions. To understand this phenomenon, we investigated the PDF of the stream-wise velocity by detailed simulations of the flow in the pore space of random sphere packs. These PDF show strongly skewed distributions with tails up to four times the average pore velocity. Negative velocities are more likely to delay tracer transport than to contribute to upstream transport of tracer.

The determination of the inhomogeneous permeability field can contribute to understand and predict the large-scale tracer dispersion. We presented an adjoint-based optimisation algorithm to estimate permeability distributions from point measurements of the velocity in a porous medium. The results show a satisfying agreement between input and estimated permeability fields. As expected, they also reveal a dependency on the observation operator.

Tracer transport on a large scale is often convection dominated. In these situations,  upstream discretisations are used which introduce additional numerical diffusivity to reduce oscillations in the solution. However, this numerical diffusivity is not always a viable solution as it strongly smears out the sharp gradients in the tracer field. In this paper, we presented a numerical analysis of the $p$-FEM method to determine under which conditions unphysical oscillations can be damped by the use of higher order methods without introducing unwanted numerical diffusion.

\bibliographystyle{unsrt}
\bibliography{paper}



\printindex
\end{document}